\documentclass[english]{article}

\usepackage{babel}
\usepackage{amstext}
\usepackage{amsmath}
\usepackage{amsfonts}
\usepackage{latexsym}
\usepackage{ifthen}
\usepackage{xypic}
\xyoption{all}
\begin{document}

%%%%%%%%%%%%%%%%%%

\newcommand{\id}{{\rm id}}
\renewcommand{\O}{{\cal O}}
\newcommand{\ball}{{\mathbb B}}
\newcommand{\Sie}{{\mathbb H}}
\newcommand{\Fam}{{\mathfrak X}}
\newcommand{\PB}{{\mathbb P}}
\newcommand{\kod}{\kappa}
\newcommand{\rk}{{\rm rk}}
\newcommand{\codim}{{\rm codim}}
\newcommand{\tr}{{\rm tr}}
\newcommand{\End}{{\cal E}nd}
\newcommand{\Hom}{{\rm Hom}}
\newcommand{\Ext}{{\rm Ext}}
\newcommand{\im}{{\rm im}}
\renewcommand{\Im}{{\mathfrak Im}}
\newcommand{\PN}[1]{{\mathbb P}_{#1}}
\newcommand{\lra}{\longrightarrow}
\newcommand{\KQ}{{\mathbb Q}}
\newcommand{\KR}{{\mathbb R}}
\newcommand{\KZ}{{\mathbb Z}}
\newcommand{\KC}{{\mathbb C}}
\newcommand{\KN}{{\mathbb N}}

\renewcommand{\theequation}{\arabic{section}.\arabic{equation}}
\newcommand{\Formel}[1]{(\ref{#1})}
\newtheorem{lemma1}[equation]{}

\newenvironment{lemma}{\begin{lemma1}{\bf Lemma.}}{\end{lemma1}}
\newenvironment{example}{\begin{lemma1}{\bf Example.}\rm}{\end{lemma1}}
\newenvironment{abs}{\begin{lemma1}\rm}{\end{lemma1}}
\newenvironment{theorem}{\begin{lemma1}{\bf Theorem.}}{\end{lemma1}}
\newenvironment{theorem2}[1]{\begin{lemma1}{\bf Theorem [#1].}}{\end{lemma1}}
\newenvironment{proposition}{\begin{lemma1}{\bf Proposition.}}{\end{lemma1}}
\newenvironment{corollary}{\begin{lemma1}{\bf Corollary.}}{\end{lemma1}}
\newenvironment{remark}{\begin{lemma1}{\bf Remark.}\rm}{\end{lemma1}}
\newenvironment{definition}{\begin{lemma1}{\bf Definition.}}{\end{lemma1}}

\newenvironment{proof}{\noindent {\em Proof}.}{\hfill $\Box$}
\newenvironment{proof2}[1]{\noindent {\em Proof of #1}.}{\hfill $\Box$}

\renewcommand{\labelenumi}{\arabic{enumi}.)}
\renewcommand{\labelenumii}{\arabic{enumi}.\arabic{enumii}.)}

\newcommand{\Lemma}[1]{Lemma~\ref{#1}}
\newcommand{\Ex}[1]{Example~\ref{#1}}
\newcommand{\Abs}[1]{\ref{#1}}
\newcommand{\Theo}[1]{Theorem~\ref{#1}}
\newcommand{\Prop}[1]{Proposition~\ref{#1}}
\newcommand{\Cor}[1]{Corollary~\ref{#1}}
\newcommand{\Rem}[1]{Remark~\ref{#1}}
\newcommand{\Def}[1]{Definition~\ref{#1}}

%%%%%%%%%%%%%%%%%

\title{On manifolds with holomorphic normal projective connections}
\author{Priska Jahnke and Ivo Radloff \thanks{The authors were supported by a Forschungsstipendium of the Deutsche For\-schungs\-gemeinschaft and the DFG--Schwerpunkt {\em Globale Methoden in der komplexen Geometrie.}}}
\date{August 2002}
\maketitle

\section*{Introduction}

Complex projective space $\PN{n}$, \'etale quotients of complex tori and compact complex manifolds whose universal cover is the unit ball ${\ball}^n \subset \KC^n$ are standard examples of complex K\"ahler manifolds admitting a (flat) holomorphic normal projective connection. In particular, any compact complex curve admits a (flat) holomorphic normal projective connection. In {\em Holomorphic projective structures on compact complex surfaces I} and {\em II}, \cite{KoOc}, Kobayashi and Ochiai proved that the list of compact complex K\"ahler surfaces admitting a normal holomorphic projective connection is precisely this list of standard examples. Their result raised the question whether or not the list is complete even in higher dimensions.

In this article we give a complete classification of projective threefolds admitting a holomorphic normal projective connection. The result shows in particular that the above list is not complete in general:

\vspace{0.2cm}

\noindent {\bf \Theo{3dim} }{\em The class of $3$--dimensional complex projective manifolds admitting a holomorphic normal projective connection consists precisely of
 \begin{enumerate}
   \item $\PN{3}$,
   \item \'etale quotients of abelian threefolds,
   \item \'etale quotients of smooth modular families of false elliptic curves,
   \item manifolds, whose universal cover is the unit ball ${\ball}^3 \subset \KC^3$.
 \end{enumerate}
\noindent However, as in the case of curves and surfaces, this list coincides with the list of projective threefolds admitting a \underline{flat} holomorphic normal projective connection.}

\vspace{0.2cm}

\noindent Recall that a false elliptic curve is an abelian surface, where the $\KQ$--endo\-morphism algebra ${\rm End}_{\KQ}$ is a totally indefinite quaternion algebra. The moduli scheme of such a surface is known to be a Shimura curve; we briefly recall the construction of the universal family in \Ex{kuga}. Modular families of false elliptic curves are well--known, but seem not to have been considered as a source of examples of manifolds admitting a flat holomorphic normal projective connection. Note that their Kodaira dimension is one. The proof of \Theo{3dim} bases on Mori theory, which is so far only sufficiently settled in the projective case, and results from variation of Hodge structures. 

In {\em On Fano manifolds with normal projective connections}, \cite{Ye}, Ye proved, that $\PN{n}$ is the only Fano manifold with a holomorphic normal projective connection. Recall that a complex (projective) manifold $X$ is called {\em Fano} if the dual of the canonical bundle $K_X$ is ample. It is called {\em minimal} if $K_X$ is nef, i.e., if $K_X$ has nonnegative intersection number with any irreducible curve in $X$. The following general structure \Theo{min} is one of the keys to the proof of \Theo{3dim}:

\vspace{0.2cm}

\noindent {\bf \Theo{min} }{\em Let $X$ be a projective manifold of dimension $n$ with a holomorphic normal projective connection. If $X \not\simeq \PN{n}$, then $X$ is minimal and does not contain any rational curve. Furthermore:
 \begin{enumerate}
  \item If $K_X \equiv 0$, then $X$ is covered by a torus.
  \item If $K_X$ is big, then $K_X$ is ample and $X$ is covered by the unit ball.
 \end{enumerate}
\noindent In general, if $K_X$ is abundant, then the Iitaka fibration $f: X \to Y$ is equidimensional and the general fiber is covered by an abelian variety. Moreover, $(n+1)^rc_r(X) = \left({n+1 \atop r}\right) c_1^r(X)$ in $H^r(X, \Omega_X^r)$ and $c_i(X) = 0$ for $i > \dim Y$.}

\vspace{0.2cm}

\noindent Recall that the {\em abundance conjecture} claims that $K_X$ is abundant if $X$ is minimal, i.e., that $|mK_X|$ is base point free for some $m \gg 0$. The induced map $f:X \to Y$ (after Stein factorization) is called {\em Iitaka fibration}. The abundance conjecture is known to hold true in dimension $\le 3$ (\cite{Ka}).

\Theo{min} is proved in the sections~\ref{sec rat} and \ref{sec min}. In section~\ref{sec surf} we give a new and shorter proof of Kobayashi and Ochiai's result of the surface case. In section~\ref{sec 3dim} the main \Theo{3dim} is proved.

\

\noindent {\bf Acknowledgements.} The idea to this paper came up during a stay of the authors at the University of Michigan, Ann Arbor, USA. The authors are grateful to Professor R. Lazarsfeld and Professor Th. Peternell for several encouraging discussions. The authors are grateful to Professor E. Viehweg for explaining one of his recent results.

\newpage

%%%%%%%%%%%%%%%%%%

\section{Projective structures and connections} \label{sec psc}

In this section we briefly recall some basic definitions and results. Throughout this section, we assume that $X$ is a compact K\"ahler manifold of dimension $n$, even though the K\"ahler condition is not always necessary. Concerning our notations see \Abs{notation}.

\begin{abs} \label{pstructure}
{\bf Holomorphic affine and projective structures.} (\cite{KoOc} or \cite{KoWu}). $X$ is said to admit a holomorphic {\em affine} (resp. {\em projective}) structure, if it can be covered by coordinate charts, such that the coordinate change is given by restrictions of holomorphic affine (resp. projective) transformations of $\KC^n$ (resp. $\PN{n}$). A manifold with a holomorphic affine structure admits a holomorphic projective structure. The ``list of standard examples'' of compact complex K\"ahler (--Einstein) manifolds admitting a holomorphic projective structure, already mentioned above, consists of
\begin{enumerate}
 \item the projective space $\PN{n}$,
 \item \'etale quotients of complex tori,
 \item manifolds, whose universal cover is the unit ball ${\ball}^n \subset \KC^n$.
\end{enumerate}
In 3. note that the group of automorphisms of ${\ball}^n$ is ${\rm SU}(1,n)$, a subgroup of ${\rm PGl}(n+1)$. By the uniformization theorem, every compact complex curve is in the list of standard examples, i.e., admits a holomorphic projective structure. The manifolds in 2. admit a holomorphic affine structure. 
\end{abs}

There is one more (projective) example in dimension $3$: the Kuga fiber space associated to a totally indefinite quaternion algebra. Even though these fiber spaces have been studied from many points of view by Kuga, Shimura, Satake, Mumford et al., they seem not to have been considered as a source of examples of manifolds admitting a holomorphic projective structure. Note that this example will not be called ``standard''. In our brief description we follow \cite{LaBi}:

\begin{example} \label{kuga}
Let $F$ be a totally indefinite quaternion algebra over $\KQ$. As an algebra, $F$ is generated by two elements $u,v \in F$ satisfying
 \[u^2 = a, \quad v^2 = b, \quad uv = -vu\]
for some $a,b \in \KQ$. $F$ is a division algebra and $F \otimes_{\KQ} \KR \simeq M_2(\KR)$. We may regard the elements of $F$ as matrices. Let $\Lambda \subset F$ be some rank $4$ $\KZ$--lattice in $F$, and let $Q \in F$ be a nondegenerated matrix satisfying $Q^t = -Q$ and $\tr(\Lambda Q \Lambda^t) \subset \KZ$. We have $Q = \left({\hspace{0.1cm} 0 \hspace{0.4cm} \alpha \atop -\alpha \hspace{0.3cm} 0}\right)$ for some $\alpha \in \KR$, $\alpha \not= 0$. 

Denote by $\Sie = \{\tau \in \KC \mid \Im(\tau) > 0\}$ the Siegel upper half plane. For any $\tau \in \Sie$, define $j_{\tau}: F \otimes_{\KQ} \KR \to {\KC}^2$ by $A \mapsto A \cdot \left({\tau \atop \alpha}\right)$. The map $j_{\tau}$ is an isomorphism of $\KR$--vector spaces for any $\tau$; it endows $F \otimes_{\KQ} \KR \simeq {\KR}^4$ with a complex structure. The lattice $\Lambda$ acts freely on ${\KC}^2 \times \Sie$ by $(z,\tau) \mapsto (z + j_{\tau}(\lambda), \tau)$. The quotient
 \[\Fam_F = ({\KC}^2 \times \Sie)/\Lambda\]
is a complex manifold and a smooth family of abelian varieties over $\Sie$. Denote the fiber of $\Fam$ over $\tau \in \Sie$ by $\Fam_{\tau}$. The lattice $\Lambda$ determines an arithmetic subgroup of the symplectic group
 \[\Gamma' = \{g \in {\rm Sp}_2(\KR) \mid \Lambda \mbox{$\left({1 \hspace{0.2cm} 0 \atop 0 \hspace{0.2cm} \alpha}\right)$} g \mbox{$\left({1 \hspace{0.2cm} 0 \atop 0 \hspace{0.2cm} \alpha}\right)^{-1}$} \subset \Lambda\}.\]
There exists a subgroup $\Gamma \subset \Gamma'$ of finite index that is torsion free (\cite{Satake}, IV, Lemma~7.2). The group $\Gamma$ acts on $\Sie$ by
 \[\tau \mapsto \frac{a\tau + b}{c\tau + d}.\]
Moreover, if $\left(a \hspace{0.2cm} b \atop c \hspace{0.2cm} d\right) \in \Gamma$, then multiplication with $(c\tau + d)$ induces an isomorphism $\Fam_{\frac{a\tau + b}{c\tau + d}} \to \Fam_{\tau}$. The induced action of $\Gamma$ on $\Fam_F \to \Sie$ is properly and discontinously. The quotient $\Fam_F/\Gamma \to \Sie/\Gamma$ is a smooth abelian fibration. Since $F$ is a division algebra, $\Sie/\Gamma$ is compact. The form $Q$ defines an embedding of $\Fam_F/\Gamma$ into projective space. If we view $\KC^2 \times \Sie$ as part of a standard coordinate chart of $\PN{3}$, then it is clear that $\Fam_F/\Gamma$ has a holomorphic projective structure. By construction, $\Fam_F/\Gamma$ is not one of the standard examples.
\end{example}

The corresponding infinitesimal notions to holomorphic affine and projective structures are holomorphic affine and projective connections. Before we come to this we recall the definition of the Atiyah class of a holomorphic vector bundle:

\begin{abs} \label{Atiyahclass}
{\bf The Atiyah class and $a(E)$.} (\cite{Ati}). Let $E$ be a holomorphic vector bundle of rank $r$ on $X$ and $\{U_{\alpha}; z_{\alpha}^1, \dots, z_{\alpha}^n\}_{\alpha \in I}$ coordinate charts with coordinates $z_{\alpha}^i$ where $E$ is trivial. Let $\{U_{\alpha}; e_{\alpha}^1, \dots, e_{\alpha}^r\}_{\alpha \in I}$ be a local frame for $E$. Let $g_{\alpha\beta} \in H^0(U_{\alpha} \cap U_{\beta}, {\rm Gl}(r, \O_X))$ be transition functions of $E$ such that $e_{\beta}^k = \sum_l g_{\alpha\beta}^{lk}e_{\alpha}^l$.

The {\em Atiyah class of $E$} is the splitting obstruction of the first jet sequence
  \[\quad 0 \lra \Omega^1_X \otimes E \lra J_1(E) \lra E \lra 0,\]
i.e., it is the image of $\id_E$ under the first connecting morphism 
 \[H^0(X, \End(E)) \lra \Ext^1(E, \Omega_X^1 \otimes E) \simeq H^1(X, \End(E) \otimes \Omega^1_X).\] 
The Dolbeault isomorphism $H^1(X, \End(E) \otimes \Omega_X^1) \simeq H^{1,1}(X, \End(E))$ maps the Atiyah class to $[-\Theta_h]$, where $\Theta_h$ denotes the canonical curvature of $E$ with respect to a hermitian metric $h$ on $E$. In particular, the trace of the Atiyah class is $-2\pi i c_1(E)$ in $H^1(X, \Omega^1_X)$.

If we define $a(E)$ as $-\frac{1}{2\pi i}$ times the Atiyah class of $E$, then the trace of $a(E)$ is $c_1(E)$, which makes this definition convenient for our purposes. In local coordinates, $a(E)$ is the class of the Chech cocycle $a(E)_{\alpha\beta} \in Z^1({\cal U}, \End(E)\otimes \Omega_X^1)$, where
 \[a(E)_{\alpha\beta} = \mbox{$\frac{1}{2\pi i} \sum_{i,j,l} \frac{\partial g_{\alpha\beta}^{jl}}{\partial z_{\alpha}^i}dz_{\alpha}^i \otimes e_{\alpha}^j \otimes e_{\beta}^{l*} = \frac{1}{2\pi i} \sum_{1 \le j,l \le r} dg_{\alpha\beta}^{jl} e_{\alpha}^j \otimes e_{\beta}^{l*}$}.\]
See \cite{Ati} for the functorial behavior of $a(E)$ under pull--back, tensor products and direct sums.
\end{abs}

\begin{abs} \label{pc}
{\bf Holomorphic affine and projective connections.} (\cite{KoOc} or \cite{KoWu}) $X$ is said to have a {\em holomorphic affine connection} if
 \[a(\Omega^1_X) = 0 \;\mbox{ in } \; H^1(X, \Omega_X^1 \otimes T_X \otimes \Omega_X^1),\]
where $T_X$ and $\Omega_X^1$ denote the holomorphic tangent sheaf and the sheaf of holomorphic $1$--forms, respectively. If $X$ has a holomorphic affine structure, then $X$ has a (flat) holomorphic affine connection. Since we assume $X$ K\"ahler, the existence of a holomorphic affine connection implies that $X$ is covered by a torus (\cite{KoWu}, 2.4.1.~Theorem). In other words: a compact K\"ahler manifold admits a holomorphic affine connection if and only if it is covered by a torus.

$X$ is said to have a {\em holomorphic (normal) projective connection} if
 \[a(\Omega_X^1) =  \mbox{$\id_{\Omega_X^1} \otimes \frac{c_1(K_X)}{n+1} + \frac{c_1(K_X)}{n+1} \otimes \id_{\Omega_X^1}$} \; \mbox{ in } \; H^1(X, \Omega_X^1 \otimes T_X \otimes \Omega_X^1),\]
where we use the identifications
 \[\End(\Omega^1_X) \otimes \Omega_X^1 \simeq \Omega_X^1 \otimes T_X \otimes \Omega_X^1 \simeq \Omega_X^1 \otimes \End(\Omega^1_X)\]
and consider $c_1(K_X)$ as an element in $H^1(X, \Omega_X^1)$. (See \cite{MoMo} for a more differential geometric description of projective connections). From now on we will drop the appellation ``normal''. If $X$ has a holomorphic projective connection and $c_1(X) = 0$ in $H^1(X, \Omega_X^1)$, then $X$ has a holomorphic affine connection. Hence: a compact K\"ahler manifold with a holomorphic projective connection and $c_1(X) = 0$ is covered by a torus.

If $X$ has a holomorphic affine (projective) structure, then $X$ admits a holomorphic affine (projective) connection. A holomorphic affine (projective) connection is said to be {\em flat} or {\em integrable} if it corresponds to a holomorphic affine (projective) structure. The examples in \Abs{pstructure} and \ref{kuga} all have a flat holomorphic projective connection. If $X$ has a holomorphic affine (projective) connection and $\tilde{X} \to X$ is \'etale, then $\tilde{X}$ admits a holomorphic affine (projective) connection. 

Gunning's formula on the Chern classes of a K\"ahler manifold with a holomorphic projective connection says (\cite{Gun}, p.94)
  \begin{equation} \label{ChernClss}
    (n+1)^rc_r(X) = \mbox{$\left({n+1 \atop r}\right)$} c_1^r(X) \; \mbox{ in } \; H^r(X, \Omega_X^r),
  \end{equation}
where $c_r(X) = c_r(T_X) \in  H^r(X, \Omega_X^r)$. Important for our purposes will be moreover the following result of Kobayashi and Ochiai:

\begin{theorem2}{Kobayashi, Ochiai} \label{KE}
 The list of K\"ahler--Einstein manifolds admitting a holomorphic projective connection is the list of standard examples.
\end{theorem2}
\Theo{KE}  is proved by showing, that $X$ is of constant holomorphic sectional curvature (see \cite{KoOc} or \cite{KoWu}, I, 1.7.1.).
\end{abs}

\begin{abs} \label{method}
{\bf Projective connections and exact sequences.} The following argument, in slightly varying form, will be used in most of the proofs: let $\nu: Y \to X$ be some morphism (not necessarily surjective), $X, Y$ compact K\"ahler manifolds, $\dim(X) = n \ge 2$. Assume that we have a bundle sequence
 \begin{equation} \label{seqKQ}
  0 \lra K \lra \nu^*\Omega_X^1 \stackrel{\rho}{\lra} Q \lra 0
 \end{equation}
on $Y$. We get maps
\[\xymatrix{
  & & H^1(Y, Q \otimes Q^* \otimes \Omega_Y^1) \ar[d]^{\mbox{$\psi$}} \\
 H^1(Y, \nu^*\Omega_X^1 \otimes \nu^*T_X \otimes \Omega_Y^1) \ar[rr]^{\mbox{$\varphi$}} & & H^1(Y, Q \otimes \nu^*T_X \otimes \Omega_Y^1),}\]
where $\varphi$ is given by $\rho \otimes \id \otimes \id$ and $\psi = \id \otimes \rho^t \otimes \id$. In this diagram, the classes $a(\nu^*\Omega_X^1)$ and $a(Q)$ are mapped onto the same class in $H^1(Y, Q \otimes \nu^*T_X \otimes \Omega_Y^1)$, i.e.,
 \begin{equation} \label{phipsi}
   \varphi(a(\nu^*\Omega_X^1)) = \psi(a(Q)).
 \end{equation}
This can be seen either by thinking of the transition functions of $\nu^*\Omega_X^1$ as upper triangular matrices corresponding to \Formel{seqKQ} or by diagram chase in the corresponding jet sequences (up to our factor $-2\pi i$, the classes in \Formel{phipsi} are the image of $\rho$ under $H^0(Y, Q \otimes \nu^*T_X) \lra H^1(Y, Q \otimes \nu^*T_X \otimes \Omega_Y^1)$).

If $X$ has a holomorphic projective connection, then we have $a(\nu^*\Omega_X^1) = \id_{\nu^*\Omega_X^1} \otimes \frac{c_1(\nu^*K_X)}{n+1} + \frac{\nu^*c_1(K_X)}{n+1} \otimes d\nu$. Here, we carefully distinguish between $c_1(\nu^*K_X)$ and $\nu^*c_1(K_X)$. Now we make two assumptions. The first is
 \[\nu^*c_1(K_X) \in \im(H^1(Y, K) \lra H^1(Y, \nu^*\Omega_X^1)).\] 
Then $\varphi(\frac{\nu^*c_1(K_X)}{n+1} \otimes d\nu) = 0$ and therefore $\varphi(a(\nu^*\Omega_X^1)) = \rho \otimes \frac{c_1(\nu^*K_X)}{n+1}$. On the other hand, $\id_Q \otimes \frac{c_1(\nu^*K_X)}{n+1} \in H^1(Y, \End(Q) \otimes \Omega_Y^1)$, and
 \[\psi(a(Q)) = \varphi(a(\nu^*\Omega_X^1)) = \mbox{$\rho \otimes \frac{c_1(\nu^*K_X)}{n+1} = \psi(\id_Q \otimes \frac{c_1(\nu^*K_X)}{n+1})$},\]
which we formally write as $\psi(a(Q(\frac{-\nu^*K_X}{n+1}))) = 0$. The second assumption is $\psi$ being injective (e.g., if \Formel{seqKQ} splits holomorphically). The injectivity implies 
 \[a(Q) = \mbox{$\id_Q \otimes \frac{c_1(\nu^*K_X)}{n+1}$} \quad \mbox{ and } \quad (n+1)^r c_r(Q) = \mbox{$\left({\rk Q \atop r}\right)$} c_1^r(\nu^*K_X).\]
Hence, if the two assumptions hold, then we get informations about the positivity of $Q$ from the positivity of $K_X$. 

Note that if $Y$ is a curve, then, by Riemann--Roch, the injectivity of $\psi$ is equivalent to the surjectivity of $H^0(Y, Q^*\otimes \nu^*\Omega_X^1) \lra H^0(Y, \End(Q)^*)$, which implies the splitting of \Formel{seqKQ}.
\end{abs}

\begin{abs} \label{notation}
\noindent {\bf Notation and conventions.} For a compact complex manifold $X$, the canonical divisor is denoted by $K_X$. If $X \to Y$ is a morphism, $K_{X/Y} = K_X - f^*K_Y$. We will identify line bundles and divisors and write $K_X$ and $K_{X/Y}$ instead of $\O_X(K_X)$ and $\O_X(K_{X/Y})$, respectively. The tensor product of line bundles will be denoted by $+$ or $\otimes$. A line bundle is called nef, if the intersection number with every irreducible curve is non--negative. It is called big, if the top self intersection class is positive. A vector bundle is called nef, if $\O_{\PB(E)}(1)$ is nef on $\PB(E)$. We write $c_i(X) = c_i(T_X)$. The Kodaira dimension of $X$ will be denoted by $\kod(X)$. The symbol $\equiv$ denotes numerical equivalence. 
\end{abs}

%%%%%%%%%%%%%%%%

\section{Manifolds containing a rational curve} \label{sec rat}
\setcounter{equation}{0}

Among the standard examples, projective space $\PN{n}$ is the only one that contains a rational curve. In this section we prove:

\begin{theorem} \label{Pn}
 Let $X$ be a compact K\"ahler manifold of dimension $n$ with a holomorphic projective connection. If $X$ contains a rational curve, then $X \simeq \PN{n}$.
\end{theorem}
Since Fano manifolds contain rational curves, we obtain Ye's result (\cite{Ye}): the only Fano manifold with a holomorphic projective connection is $\PN{n}$. Ye uses deformation theory, while our proof is completely different.

\vspace{0.2cm}

\begin{proof}
 Let $C$ be any rational curve in $X$ and $f: \PN{1} \to X$ its normalization. The claim is that $f^*T_X$ is ample. Then $X \simeq \PN{n}$ by  Mori's Theorem (\cite{Mor} and \cite{MiPe}, p.41, 4.2.~Theorem in particular). On $\PN{1}$, the bundle $f^*T_X$ splits
 \[f^*T_X = \oplus_{i=1}^n \O_{\PN{1}}(a_i), \quad a_1 \le \cdots \le a_n, \quad a_n \ge 2.\]
If all $a_i \ge 2$ then the claim is obviously true, hence assume that there exists an index $1 \le i_0 \le n-1$, such that $a_i \le 1$ for $i \le i_0$. We have the two non--trivial subbundles $T_{\le 1} = \oplus_{a_i \le 1} \O_{\PN{1}}(a_i)$ and $T_{>1} = \oplus_{a_i>1} \O_{\PN{1}}(a_i)$, and $f^*T_X = T_{\le 1} \oplus T_{>1}$. It is enough to show that $T_{\le 1}$ is ample. To this end, we apply the method from \Abs{method} to
 \begin{equation} \label{Seq Om}
  0 \lra T^*_{>1} \lra f^*\Omega_X^1 \stackrel{\rho}{\lra} T^*_{\le 1} \lra 0.
 \end{equation}
Since $H^1(\PN{1}, T^*_{\le 1}) = 0$, the first assumption in \Abs{method} is satisfied, claiming that $f^*c_1(K_X) \in \im(H^1(\PN{1}, T^*_{>1}) \lra H^1(\PN{1}, f^*\Omega_X^1))$. Hence $\psi(a(T^*_{\le 1}(\frac{-f^*K_X}{n+1}))) = 0$ by \Abs{method}. Since \Formel{Seq Om} splits, $\psi$ is injective, i.e., the second assumption in \Abs{method} is satisfied as well. Hence
 \[a(T^*_{\le 1}) = \mbox{$\id_{T_{\le 1}^*} \otimes \frac{c_1(f^*K_X)}{n+1}$} \; \mbox{ in } \; H^1(\PN{1}, \End(T^*_{\le 1}) \otimes \Omega_{\PN{1}}^1).\]
Since $T_{\ge 1}$ splits into a sum of line bundles, contraction of $a(T_{\le 1})$ to the $i$--th diagonal entry gives $c_1(\O_{\PN{1}}(a_i)) = -\frac{c_1(f^*K_X)}{n+1}$ in $H^1(\PN{1}, \Omega_{\PN{1}}^1)$ for $i = 1, \dots, i_0$. Hence, if $-K_X.C = r$, then 
 \[(n+1) \mid r \quad \mbox{ and } \quad T_{\le 1} \simeq \oplus_{i_0}\O_{\PN{1}}(\mbox{$\frac{r}{n+1}$}).\] 
This shows $\deg T_{\le 1} = i_0 \cdot \frac{r}{n+1}$. Since $\deg(T_{> 1} \oplus T_{\le 1}) = \deg f^*T_X = -K_X.C = r$, we have $\deg T_{>1} = (n+1-i_0) \cdot \frac{r}{n+1}$. But $\deg T_{> 1} > 0$ by definition of $T_{> 1}$. Hence $-K_X.C = r > 0$, and $T_{\le 1}$ is ample.
\end{proof}

%%%%%%%%%%%%

\section{Minimal manifolds} \label{sec min}
\setcounter{equation}{0}

Let $X$ be a projective manifold with a holomorphic projective connection different from $\PN{n}$. By \Theo{Pn}, $X$ does not contain a rational curve. The {\em cone theorem}, first proved by Mori (\cite{Mo82}), states that if $K_X$ is not nef, then $X$ contains a rational curve. Hence, $K_X$ is nef and $X$ is minimal. The abundance conjecture claims that moreover $|mK_X|$ is base point free for some $m \gg 1$. The induced map
 \[f: X \lra Y,\]
after Stein factorization, is called {\em Iitaka fibration}. Here, $Y$ is a normal projective variety of dimension $\kod(X)$. The abundance conjecture is known to hold true in the case $\dim X \le 3$.

\begin{theorem} \label{min}
Let $X$ be a projective manifold of dimension $n$ with a holomorphic normal projective connection. If $X \not\simeq \PN{n}$, then $X$ is minimal and does not contain any rational curve. Furthermore:
 \begin{enumerate}
  \item If $K_X \equiv 0$, then $X$ is covered by a torus.
  \item If $K_X$ is big, then $K_X$ is ample and $X$ is covered by the unit ball.
 \end{enumerate}
\noindent In general, if $K_X$ is abundant, then the Iitaka fibration $f: X \to Y$ is equidimensional and the general fiber is covered by an abelian variety. Moreover, $(n+1)^rc_r(X) = \left({n+1 \atop r}\right) c_1^r(X)$ in $H^r(X, \Omega_X^r)$ and $c_i(X) = 0$ for $i > \dim Y$.
\end{theorem}

\begin{proof}
If $X \not\simeq \PN{n}$, then $X$ does not contain a rational curve by \Theo{Pn}. By the cone theorem, $K_X$ is nef. If $K_X \equiv 0$, then $X$ is covered by a torus (see \Abs{pc}). 

Assume that $K_X$ is abundant and let $f:X \to Y$ be the Iitaka fibration. Since $mK_X$ is trivial on each fiber of $f$, the dual of $K_X$ is $f$--nef. Hence, by a theorem of Mori and Mukai (\cite{MoMu}, Theorem~2), any irreducible component of ${\rm Exc}(f) = \{x \in X \mid \dim_x f^{-1}(f(x)) > \dim X - \dim Y\}$ is covered by rational curves, if it is not empty. Since $X$ does not contain a rational curve by \Theo{Pn}, ${\rm Exc}(f) = \emptyset$, i.e., $f$ is equidimensional. 

If $K_X$ is big, then $K_X$ is abundant by the base point free theorem (see for example \cite{KMM}). The Iitaka fibration $f$ is hence a birational, equidimensional morphism. Hence $K_X$ is ample and $f$ is an isomorphism. By the theorem of Aubin and Yau, $X$ admits a K\"ahler--Einstein metric. Claim 2. hence follows from \Theo{KE}. Claim 1. and 2. are proved.

Since $\dim Y = \kod(X) = \nu(X) = \max_{\nu}\{c_1^{\nu}(K_X) \not= 0\}$, the Chern classes $c_i(X)$ vanish for $i > \dim Y$ by Gunning's formula \Formel{ChernClss}. Let $F$ be the general fiber of $f$. Then $N_{F/X} \simeq \O_F^{\oplus \dim Y}$, and the adjunction formula shows $c_1(K_F) = 0$. The following \Lemma{genfiber} completes the proof of \Theo{min}.
\end{proof}

\begin{lemma} \label{genfiber}
 Let $X$ be a projective manifold of dimension $n$ with a holomorphic projective connection. Let $f: X \to Y$ be some morphism with connected fibers, $Y$ normal. If the canonical bundle of the general fiber $F$ is numerically trivial, then $F$ is covered by an abelian variety.
\end{lemma}

\begin{proof}
There exists a finite unramified cover $\tilde{\nu}: A \times N \to F$, where $A$ is abelian and $N$ is simply connected with vanishing first Chern class, \cite{Beau}. We have to show $\dim N = 0$. Assume $\dim N = r > 0$. Identify $N$ and $\{a\} \times N$ for some fixed $a \in A$ and let $\nu: N \to F$ be the induced map. We have the exact sequence
 \begin{equation} \label{omn}
   0 \lra \O_N^{\oplus (n-r)} \lra \nu^*\Omega_X^1 \lra \Omega_N^1 \lra 0
 \end{equation}
on $N$. By definition, $a(\nu^*\Omega_X^1) = \id_{\nu^*\Omega_X^1} \otimes \frac{c_1(\nu^*K_X)}{n+1} + \frac{\nu^*c_1(K_X)}{n+1} \otimes d\nu$, since $X$ has a holomorphic projective connection. Using the maps from \Abs{method} defined by \Formel{omn}, we find $\varphi(a(\nu^*\Omega_X^1)) = d\nu \otimes \frac{c_1(\nu^*K_X)}{n+1} + \frac{c_1(\nu^*K_X)}{n+1} \otimes d\nu$. Define
 \[\xi = \mbox{$\id_{\Omega_N^1} \otimes \frac{c_1(\nu^*K_X)}{n+1} + \frac{c_1(\nu^*K_X)}{n+1} \otimes \id_{\Omega_N^1}$} \in H^1(N, \Omega_N^1 \otimes T_N \otimes \Omega_N^1).\]
Then $\psi(\xi) = d\nu \otimes \frac{c_1(\nu^*K_X)}{n+1} + \frac{c_1(\nu^*K_X)}{n+1} \otimes d\nu$, which implies
 \[\psi(a(\Omega_N^1)) = \varphi(a(\nu^*\Omega_X^1)) = \psi(\xi).\]
If we assume, as in \Abs{method}, that $\psi$ is injective, then
 \begin{equation} \label{atomn}
  a(\Omega_N^1) = \xi = \mbox{$\id_{\Omega_N^1} \otimes \frac{c_1(\nu^*K_X)}{n+1} + \frac{c_1(\nu^*K_X)}{n+1} \otimes \id_{\Omega_N^1}$}.
 \end{equation}
The trace gives $\tr(a(\Omega_N^1)) = c_1(K_N) = \frac{r}{n+1}c_1(\nu^*K_X) + \frac{1}{n+1}c_1(\nu^*K_X)$. Hence $\frac{c_1(K_N)}{r+1} = \frac{c_1(\nu^*K_X)}{n+1}$, which means by \Formel{atomn}, that $N$ has a holomorphic projective connection. Since $c_1(K_N) = 0$, $N$ is covered by a torus by \Abs{pc}, contradicting the fact, that $N$ is simply connected. It remains to show the injectivity of $\psi$, which is equivalent to the surjectivity of 
 \[H^0(N, \Omega_N^1 \otimes \nu^*T_X \otimes \Omega_N^1) \lra H^0(N, \Omega_N^1 \otimes \O_N^{\oplus (n-r)} \otimes \Omega_N^1).\] 
All elements in question are symmetric, i.e., $a(\Omega_N^1), \xi \in H^1(N, S^2\Omega_N^1 \otimes T_N)$. It hence suffices to show the surjectivity of 
 \[H^0(N, S^2\Omega_N^1 \otimes \nu^*T_X) \lra H^0(N, S^2\Omega_N^1 \otimes \O_N^{\oplus (n-r)}).\]
But $H^0(N, S^2\Omega_N^1) = 0$ by \cite{Peternell}, Theorem~5.6. and \cite{Tian}, Theorem~2.1. The lemma is proved.
\end{proof}

\

We conclude this section by applying method \Abs{method} to submanifolds with splitting tangent sequence. A theorem of Van de Ven states, the only compact submanifolds of $\PN{n}$ with holomorphically splitting tangent sequence are linear subspaces, \cite{VdV}.

\begin{proposition} \label{split}
Let $X$ be an $n$--dimensional compact K\"ahler manifold with a holomorphic projective connection. If $Y$ is an $m$--dimensional compact submanifold with splitting tangent sequence, then $Y$ admits a holomorphic projective connection and 
 \begin{equation} \label{c1XY}
  \mbox{$\frac{c_1(K_Y)}{m+1} = \frac{c_1(K_X|_Y)}{n+1}$} \; \mbox{ in } \; H^1(Y, \Omega_Y^1).
 \end{equation}
\end{proposition}

\begin{proof}
By assumption, the canonical sequence
 \begin{equation} \label{xysplit}
  0 \lra N^*_{Y/X} \lra \Omega_X^1|_Y \lra \Omega_Y^1 \lra 0
 \end{equation}
splits holomorphically. As in the proof of \Lemma{genfiber} we get $a(\Omega^1_Y) = \id_{\Omega^1_Y} \otimes \frac{c_1(K_X\mid_Y)}{n+1} + \frac{c_1(K_X\mid_Y)}{n+1}\otimes \id_{\Omega^1_Y}$. The trace gives \Formel{c1XY}; hence $Y$ has a holomorphic projective connection.
\end{proof}

\

Note that Van de Ven's theorem follows immediately: let $Y$ be a submanifold of $\PN{n}$ with splitting tangent sequence. By \Formel{c1XY} from \Prop{split}, $-K_Y$ is ample, i.e., $Y$ is Fano with a holomorphic projective connection. By \Theo{Pn}, $Y \simeq \PN{m}$. If $d = \deg(Y)$, then \Formel{c1XY} implies $d = 1$, i.e., $Y$ is linearly embedded (note, however, that it is not difficult to conclude Van de Ven's theorem directly from Mori's theorem on manifolds with ample tangent bundle).

%%%%%%%%%%%%%%

\section{Surfaces} \label{sec surf}
\setcounter{equation}{0}

In \cite{KoOc}, Kobayashi and Ochiai give a list of all classes of compact complex surfaces admitting a holomorphic projective connection. In this section, we give a short alternative proof of this result in the case of K\"ahler surfaces since it shows the general principle of how to deal with minimal manifolds admitting a holomorphic projective connection in case of abundance (\Prop{ellsurf}).

\begin{theorem} \label{surf}
  Every compact K\"ahler surface with a holomorphic projective connection is in the list of standard examples. 
\end{theorem}

\begin{remark}
In the non--K\"ahler case, certain Hopf surfaces and other surfaces with a holomorphic affine connection add to the list (\cite{KoOc} and \cite{IKO}).
\end{remark}

\noindent {\em Proof of \Theo{surf}.} Let $X$ be a compact complex K\"ahler surface with a holomorphic projective connection, $X$ not necessarily projective. By \Theo{Pn}, $X$ contains no $(-1)$--curve, i.e., $X$ is minimal. Assume that $X$ is not a standard example. By \Formel{ChernClss}, $c_1^2 = 3c_2$. The list of minimal K\"ahler surfaces (see for example \cite{BPV}) shows that the only case possibly remaining is that of a properly elliptic surface where $\kod(X) = 1$. \Theo{surf} is therefore proved by:

\begin{proposition} \label{ellsurf}
  A minimal properly elliptic K\"ahler surface does not admit a holomorphic projective connection.
\end{proposition}

\begin{proof}
Let $X$ be a minimal properly elliptic K\"ahler surface and let $f:X \to C$ be the morphism defined by $|mK_X|$ for some $m \gg 0$, $C$ a smooth compact curve. Assume $X$ admits a holomorphic projective connection. By \Formel{ChernClss}, $c_1^2 = c_2 = 0$. By \Theo{Pn}, the only singular fibers of $f$ are of type ${}_mI_0$, i.e., multiples of a smooth fiber. After an \'etale cover we may assume $f_*K_{X/C} \simeq \O_C$ (\cite{BPV}, III, 18.2. and 18.3.). Define $D = \sum_{i=1}^l (m_i-1)F_i$, where $F_i$ are the singular fibers with multiplicities $m_i$, $i=1, \dots, l$. The canonical bundle formula (\cite{BPV}, V, 12.1.) gives
 \[K_X \simeq f^*K_C \otimes \O_X(D).\]
The differential $df:f^*K_C \to \Omega_X^1$ factorizes over $f^*K_C\otimes \O_X(D)$ and we obtain the exact sequence
 \begin{equation} \label{tangseq}
  0 \lra f^*K_C \otimes \O_X(D) \lra \Omega_X^1 \stackrel{\rho}{\lra} \O_X \lra 0.
 \end{equation} 
From $\chi(\O_X) = \frac{1}{24}(c_1^2 + 2c_2) = 0$ we infer $q(X) = h^0(X, K_X) + 1$. Hence \Formel{tangseq} is exact on $H^0$--level, implying that the sequence splits holomorphically. Then, by \Abs{method} using sequence \Formel{tangseq}, one has $c_1(K_X) = 3 c_1(\O_X) = 0$, contradicting $\kod(X) = 1$. The proofs of \Prop{ellsurf} and \Theo{surf} are complete.
\end{proof}

%%%%%%%%%%%%%%%%%

\section{Threefolds} \label{sec 3dim}
\setcounter{equation}{0}

Let $X$ be a smooth projective variety of dimension $3$ with holomorphic projective connection, different from $\PN{3}$. By \Theo{min}, $X$ is minimal and $X$ is in the list of standard examples, provided $\kod(X) \not = 1, 2$. Since $K_X$ is nef and the abundance conjecture holds in dimension $3$, we have the Iitaka fibration 
 \[f: X \lra Y,\]
induced by $|mK_X|$ for some $m \gg 1$. By \Theo{min}, $f$ is equidimensional and the general fiber is covered by an abelian variety. Since $X$ does not contain a rational curve by \Theo{Pn}, $f$ has a very special structure:
\begin{itemize}
  \item If $\kod(X) = 1$, then $f$ is an almost smooth abelian or hyperelliptic fibration over a smooth curve $Y = C$, i.e., the only singular fibers are multiples of an abelian or hyperelliptic surface (\cite{Oguiso}, Theorem~B.1.). 
  \item If $\kod(X) = 2$, then $f$ is an almost smooth elliptic fibration over a normal surface $Y = S$, i.e., the only singular fibers are multiples of an elliptic curve.
\end{itemize}
The aim of this section is to prove the main theorem:

\begin{theorem} \label{3dim}
The class of $3$--dimensional complex projective manifolds admitting a holomorphic normal projective connection consists exactly of
 \begin{enumerate}
   \item $\PN{3}$,
   \item \'etale quotients of abelian threefolds,
   \item \'etale quotients of smooth modular families of false elliptic curves,
   \item manifolds, whose universal cover is the unit ball ${\ball}^3 \subset \KC^3$.
 \end{enumerate}
\noindent This list coincides with the list of projective threefolds admitting a flat holomorphic normal projective connection.
\end{theorem}

\begin{proof}
The examples \Abs{pstructure} and \ref{kuga} in section~\ref{sec psc} prove that any manifold in \Theo{3dim} has a (flat) holomorphic projective connection. Conversely, let $X$ be a projective threefold with a holomorphic projective connection. If $X$ is not in the list of standard examples, then $X$ is minimal and the Iitaka fibration maps onto a normal curve or surface. The surface case $\kod(X) = 2$ is excluded by \Prop{kod2} below. In the case $\kod(X) = 1$, $X$ is, up to an \'etale cover, a smooth modular family of false elliptic curves by \Cor{fec}. 
\end{proof}

\begin{proposition} \label{kod2}
A projective threefold of Kodaira dimension $2$ does not admit a holomorphic projective connection.
\end{proposition}

\begin{proof2}{\Prop{kod2}}
Assume that $X$ is a projective threefold with $\kod(X) = 2$ admitting a holomorphic projective connection. By \Theo{min}, $X$ is minimal. Let $f: X \to S$ be the Iitaka fibration, $S$ a normal surface. $f$ is an almost smooth elliptic fibration. Let $U \subset S$ be a Zariski open set in $S$, such that $U^C$ contains $S_{sing}$ and the points in $S$, where $sing(f)$ is not a normal crossing divisor.

The double dual of some tensor power of $f_*K_{X/S}$ (a priori only defined on $S_{reg}$) is trivial and induces a cover $\tilde{S} \to S$ from a normal surface $\tilde{S}$, unramified over $U$. Let $S'$ and $X'$ be desingularizations of $\tilde{S}$ and of the fiber product $X \times_S S'$, respectively, such that $S' \to S$ and $X' \to X$ are \'etale over $U$ and $f^{-1}(U)$, respectively. Let $f': X' \to S'$ be the induced map. By construction, $f'_*K_{X'/S'} \simeq \O_{S'}$ on some Zariski open subset $U' \subset S'$, where the complement consists of points. Let $V' \subset U'$ be the Zariski open subset of $U'$, where $f'$ is smooth. The (flat) section of $f'_*K_{X'/S'}|_{V'}$ induces a section of $R^1f'_*\KC|_{V'}$. By \cite{Fuj}, Lemma~4.4.,
 \[H^0(S', R^1f'_*\KC) \simeq H^0(V', R^1f'_{V'*}\KC).\]
Leray spectral sequence gives an exact sequence
 \[0 \lra H^1(S', \KC) \lra H^1(X', \KC) \lra H^0(S', R^1f'_*\KC) \lra 0.\]
Hence $q(X') = q(S') + 1$, i.e., there exists a holomorphic $1$--form $\omega$ on $X'$, which does not vanish on the general fiber of $f'$.

Let $H' \subset S'$ be the preimage of some sufficiently general curve $H \subset S$ contained in $U$, and let $X'_{H'} = {f'}^{-1}(H')$. Then $f'_{H'}: X'_{H'} \to H'$ is an elliptic surface and $f'_{H'*}K_{X'_{H'}/H'} \simeq \O_{H'}$. The differential $\Omega^1_{X'}|_{X'_{H'}} \to \Omega^1_{X'_{H'}}$ followed by the surjection $\Omega^1_{X'_{H'}} \to \O_{X'_{H'}}$ from \Formel{tangseq} in the proof of \Prop{ellsurf} induces a sequence
 \begin{equation} \label{rhotil}
  0 \lra K' \lra \Omega^1_{X'}|_{X'_{H'}} \lra \O_{X'_{H'}} \lra 0
 \end{equation}
(if $X'$ is smooth over $U'$, then $K' = {f'}^*\Omega^1_{U'}|_{X'_{H'}})$. The $1$--form $\omega$ guaranties that \Formel{rhotil} is surjective on $H^0$--level. Sequence \Formel{rhotil} hence splits. Since $X' \to X$ is \'etale in a neighborhood of $X'_{H'}$, the bundle $\Omega^1_{X'}|_{X'_{H'}}$ is the pull--back of $\Omega^1_X$ to $X'_{H'}$. The argument from \Abs{method} shows now the pull--back of $K_X$ to $X'_{H'}$ is numerically trivial, contradicting $\kod(X) = 2$.
\end{proof2}

\

It remains to consider the case $\kod(X) = 1$. By \Ex{kuga}, a smooth projective threefold with a holomorphic projective connection and $\kod(X) = 1$ does exist. Our aim is to prove that this is the only example. Before we come to this, we briefly recall:

\begin{abs} \label{ram}
{\bf Ramified resolution of multiple fibers.} Let $f: X \to C$ be an almost smooth fibration with multiple fibers $F_1, \dots, F_l$ of multiplicities $m_1, \dots, m_l$. Define
\[D = \mbox{$\sum\nolimits_{i=1}^l$} (m_i-1)F_i.\] 
A ramified base--change $\pi: C' \to C$ of degree $m = {\rm l.c.m.}(m_1, \dots, m_l)$, ramified over the critical values of $f$ and one additional point $a_0 \in C$, leads to a smooth fibration $f':X' \to C'$, where $X'$ is the normalization of $X \times_C C'$ (see for example \cite{Kod}). Let $\mu: X' \to X$ be the induced ramified cover. On $X'$ we have the commutative diagram
 \begin{equation} \label{diag}
  \xymatrix{0 \ar[r] & \mu^*(f^*K_C \otimes \O_X(D)) \ar[r] \ar@{^{(}->}[d] & \mu^*\Omega_X^1 \ar[r]^{\mbox{$\rho$}} \ar@{^{(}->}[d] & \Omega_{X'/C'}^1 \ar[r] \ar@{=}[d] & 0\\
 0 \ar[r] & {f'}^*K_{C'} \ar[r] & \Omega_{X'}^1 \ar[r] & \Omega_{X'/C'}^1 \ar[r] & 0,}
 \end{equation}
where $\rho$ is the pull--back of the cokernel map of $f^*K_C \otimes \O_X(D) \hookrightarrow \Omega_X^1$.
\end{abs} 

\begin{proposition} \label{kod1}
Let $X$ be a projective threefold with a holomorphic projective connection and $\kod(X) = 1$. Then $X$ is, up to an \'etale cover, a smooth abelian fibration over a smooth curve $C$ and 
 \begin{equation} \label{arakelov}
   c_1(K_X)= 2c_1(K_{X/C}) \; \mbox{ in } \; H^1(X, \Omega_X^1).
 \end{equation}
\end{proposition}

A smooth family $f: X \to C$ of abelian surfaces as in \Prop{kod1} is non--isotrivial. Let $E^{1,0} = f_*\Omega_{X/C}^1$ and $E^{0,1} = R^1f_*\O_X$. The Higgs field of the Higgs bundle $(E^{1,0} \oplus E^{0,1}, \theta)$ is given by the edge morphism
 \begin{equation} \label{theta10}
   \theta_{1,0}: f_*\Omega_{X/C}^1 \lra R^1f_*f^*K_C = K_C \otimes R^1f_*\O_X
 \end{equation}
of the sequence
 \[0 \lra f^*K_C \lra \Omega_X^1 \lra \Omega_{X/C}^1 \lra 0.\]
The equality $c_1(K_{X/C}) = c_1(f^*K_C)$ implies that $\theta_{1,0}$ is an isomorphism (\cite{ViZu} or see the proof of Claim~2), and that the Arakelov inequality (\cite{Peters}, \cite{ViZu}) is sharp:
 \[2 \deg(E^{1,0}) = 2 \deg(K_C).\]

The authors are grateful to E. Viehweg for pointing out one of his and K. Zuo's recent results (\cite{ViZu}): let $f: X \to C$ be a smooth non--isotrivial family of abelian surfaces with maximal Higgs field. Then there exists an \'etale cover $C' \to C$, such that $C'$ is a Shimura curve parametrizing false elliptic curves and $X' = X \times_C C'$ is the corresponding universal family (\cite{ViZu}, Corollary~5.2). One idea is to study the global endomorphisms of $R^1f_*\KQ$ in order to get informations on ${\rm End}_{\KQ}(A)$, $A$ the general fiber of $f$.

\Prop{kod1} hence implies

\begin{corollary} \label{fec}
 Let $X$ be a projective threefold with a holomorphic projective connection and $\kod(X) = 1$. Then $X$ is, up to an \'etale cover, a smooth modular family of false elliptic curves.
\end{corollary}
\Cor{fec} completes the proof of \Theo{3dim}.

\

\begin{proof2}{\Prop{kod1}}
Let $X$ be as in \Prop{kod1}. By \Theo{min}, $X$ is minimal. Let $f:X \to C$ be the Iitaka fibration of $X$. As noted above, $f$ is an almost smooth abelian or hyperelliptic fibration. 

The proof is in two steps. We start with a ramified resolution of the singular fibers $f': X' \to C'$ as explained in \Abs{ram} and prove a numerical property in Claim~1. Then we show that there is an \'etale resolution of the singular fibers in Claim~2 to 4. This will conclude the proof.

\vspace{0.2cm}

\noindent {\bf Claim~1:} {\em The map $f'$ is a smooth family of abelian surfaces with moduli in the fibers. Moreover,
\begin{equation} \label{claim1}
  a(\Omega_{X'/C'}^1) = \mbox{$\id_{\Omega_{X'/C'}^1} \otimes \frac{c_1(\mu^*K_X)}{4}$} \mbox{ in } H^1(X', \End(\Omega_{X'/C'}^1) \otimes \Omega_{X'}^1).
 \end{equation}} 

\noindent {\em Proof of Claim~1.} For simplicity, we first assume $f$ is smooth, i.e., $D = 0$. Note that in this case \Formel{claim1} implies the Chern class equality \Formel{arakelov}. 

Consider the canonical sequence
 \begin{equation} \label{smooth}
  0 \lra f^*K_C \lra \Omega_X^1 \lra \Omega_{X/C}^1 \lra 0.
 \end{equation}
We want to prove $a(\Omega_{X/C}^1(\frac{-K_X}{4})) = a(\Omega_{X/C}^1) - \id_{\Omega_{X/C}^1} \otimes \frac{c_1(K_X)}{4} = 0$ and that $f$ is a non--trivial family of abelian surfaces. The method from \Abs{method} applied to sequence \Formel{smooth} gives a map
 \begin{equation} \label{psi}
  \psi: H^1(X, \End(\Omega_{X/C}^1) \otimes \Omega_X^1) \lra H^1(X, \Omega_{X/C}^1 \otimes T_X \otimes \Omega_X^1),
 \end{equation}
such that $\psi(a(\Omega_{X/C}^1(\frac{-K_X}{4}))) = 0$. It is easy to see that $\psi$ is injective, if $f$ is an analytic bundle of abelian surfaces. This implies $a(\Omega_{X/C}^1(\frac{-K_X}{4})) = 0$ and hence $0 = 2c_1(K_{X/C}) = c_1(K_X)$, contradicting $\kod(X) = 1$. If the general fiber of $f$ is not abelian, then it is hyperelliptic. In this case, $X$ is an elliptic bundle over the relative Albanese variety $A(X/C)$. Here, $A(X/C)$ is projective and an elliptic bundle over $C$; the fibers are the Albanese tori of the fibers of $X \to C$ (for details see \cite{DPS}, 3.12., or \cite{Cam}). As above, we get a contradiction to $\kod(X) = 1$. 

The map $f$ hence is a smooth family of abelian surfaces with moduli in the fibers. The first part of Claim~1 is proved. Unfortunately, in this case, $\psi$ will not be injective, i.e., we need a different argument to conclude formula \Formel{claim1}. 

Since $f$ is an abelian fibration, $E = E^{1,0} = f_*\Omega_{X/C}^1$ is a nef rank $2$ vector bundle on $C$ (\cite{Griff}). Moreover, $\Omega_{X/C}^1 = f^*E$. The Kodaira--Spencer map of $f$ is non--zero, implying that the tangent sequence of a fiber is non--split in general. This shows $f_*T_X \simeq f_*T_{X/C} \simeq E^*$ and 
 \begin{equation} \label{rank}
   \rk(f_*\Omega_X^1) = 1,2.
 \end{equation}

For some $m \gg 1$, we have $mK_X = f^*A$, $A$ some ample divisor on $C$. The class $a(\Omega_{X/C}^1(\frac{-K_X}{4}))$ is the pull--back followed by the differential $df$ of $a(E(\frac{-A}{4m})) \in H^1(C, \End(E) \otimes K_C)$. The map $\psi$ from \Formel{psi} fits into the diagram obtained from Leray spectral sequence
 \[\xymatrix{H^1(C, \End(E) \otimes f_*\Omega_X^1) \ar[d]^{\mbox{$\psi'$}} \ar@{^{(}->}[r] & H^1(X', \End(\Omega_{X/C}^1) \otimes \Omega_X^1) \ar[d]^{\mbox{$\psi$}}\\
H^1(C, E \otimes f_*(T_X \otimes \Omega_X^1)) \ar@{^{(}->}[r] & H^1(X, \Omega_{X/C}^1 \otimes T_X \otimes \Omega_X^1).}\]
The top line maps the image of $a(E(\frac{-A}{4m}))$ in $H^1(C, \End(E) \otimes f_*\Omega_X^1)$ onto the class $a(\Omega_{X/C}^1(\frac{-K_X}{4}))$. Since $\psi(a(\Omega_{X/C}^1(\frac{-K_X}{4}))) = 0$, the image of $a(E(\frac{-A}{4m}))$ in $H^1(C, \End(E) \otimes f_*\Omega_X^1)$ is mapped to zero by $\psi'$. Hence, if $\psi'$ is injective, then $a(\Omega_{X/C}^1(\frac{-K_X}{4})) = 0$ as claimed. 

The map $\psi'$ is injective, if the inclusion $E^* \otimes f_*\Omega_X^1 \hookrightarrow f_*(T_X \otimes \Omega_X^1)$, obtained from the push--forward of 
 \[0 \lra T_{X/C} \otimes \Omega_X^1 \lra T_X \otimes \Omega_X^1 \lra f^*T_C \otimes \Omega_X^1 \lra 0,\]
splits holomorphically. By \Formel{rank}, $\rk(f_*\Omega_X^1) = 1,2$. A direct computation shows $\rk(f_*(T_X \otimes \Omega_X^1)) = 2 \cdot \rk(f_*\Omega_X^1) + 1$. We claim that $E^* \otimes f_*\Omega_X^1 \hookrightarrow f_*(T_X \otimes \Omega_X^1)$ defines the extension
 \begin{equation} \label{ext}
  0 \lra E^* \otimes f_*\Omega_X^1 \lra f_*(T_X \otimes \Omega_X^1) \lra \O_C \lra 0.
 \end{equation}
Indeed, first note that the cokernel of $E^* \otimes f_*\Omega_X^1 \hookrightarrow f_*(T_X \otimes \Omega_X^1)$ is contained in $T_C \otimes f_*\Omega_X^1$ and hence free. Secondly, note that the cokernel has a section coming from $\id \in H^0(X, T_X \otimes \Omega_X^1)$. Using this and the fact that $f$ is smooth, we get \Formel{ext}.

In \Formel{ext}, the canonical section in $H^0(C, f_*(T_X \otimes \Omega_X^1)) = H^0(X, T_X \otimes \Omega_X^1)$ maps to $1 \in H^0(C, \O_C)$, i.e., \Formel{ext} is surjective on $H^0$--level. This shows that \Formel{ext} splits holomorphically. Hence $\psi'$ is injective and Claim~1 is proved in the case where $f$ is smooth.

If $f$ is not smooth, i.e., if $D \not= 0$, then the proof of Claim~1 gets slightly more complicated. One has to use the top line of \Formel{diag} instead of sequence \Formel{smooth} and has to make use of the fact that $\Omega_{X'/C'}^1$ comes from $X$. However, since the main argumentation remains unchanged, we omit the details. Claim~1 is proved.

\ 

In order to complete the proof of \Prop{kod1}, it remains to show that we can find an {\em \'etale} cover of $X$, resolving the multiple fibers (which is only a problem if $C = \PN{1}$).

By Claim~1, $f'$ is a smooth fibration of abelian surfaces. As in the smooth case define $E = f'_*\Omega_{X'/C'}^1$. This is a nef rank $2$ vector bundle and
 \[\Omega_{X'/C'}^1 = {f'}^*E.\]
Let $R_D \subset C'$ be the part of the ramification divisor of $C' \to C$ lying over the critical values of $f$.

\vspace{0.2cm}

\noindent {\bf Claim~2:} {\em The edge morphism from \Formel{diag} 
 \begin{equation} \label{thetad}
   \theta_D: f'_*\Omega_{X'/C'}^1 = E \lra f'_*\mu^*(f^*K_C \otimes \O_X(D)) \otimes R^1f'_*\O_{X'}
 \end{equation}
is an isomorphism.}

\vspace{0.2cm}

\noindent If $f: X \to C$ is a smooth abelian fibration, then Claim~1 implies, that $f$ reaches the Arakelov bound. It is well--known that $\theta_D = \theta_{1,0}$ from \Formel{theta10} is an isomorphism in this case (see \cite{ViZu}). 

\vspace{0.2cm}

\noindent {\em Proof of Claim~2.} By \Formel{rank}, $\rk(f'_*\Omega_{X'}^1) = \rk(f'_*\mu^*\Omega_X^1) = 1,2$. In fact, the rank is one: with notation from \cite{ViZu}, $E = E^{1,0} = F^{1,0} \oplus N^{1,0}$, where $N^{1,0} = \ker(\theta_{1,0})$ is a numerically trivial line bundle or zero, and $F^{1,0}$ is ample. But if $E$ is decomposable, it will decompose in a sum of two line bundles of the same degree by Claim~1. Therefore, $N^{1,0} = 0$ and $\rk(f'_*\mu^*\Omega_X^1) = 1$. 

Since $\mu^*(f^*K_C \otimes \O_X(D)) = {f'}^*(\sigma^*K_C \otimes \O_{C'}(R_D))$, the top sequence in \Formel{diag} shows that $\theta_D$ is injective. Using relative duality and $c_1(K_X) = 2c_1(f^*K_C \otimes \O_X(D))$ from Claim~1, a direct computation shows immediately that the degree of the vector bundle on the right side in \Formel{thetad} is $\deg E$. The map $\theta_D$ is hence an isomorphism. Claim~2 is proved.

\

\noindent Claim~1 implies $c_1(K_X) = 2c_1(f^*K_C \otimes \O_X(D))$. We may assume that this does not only hold numerically:

\vspace{0.2cm}

\noindent {\bf Claim~3:} {\em We may assume $K_X \simeq (f^*K_C \otimes \O_X(D))^2$.}

\vspace{0.2cm}

\noindent {\em Proof of Claim~3.} The isomorphisms $\theta_D$ and $R^1f'_*\O_{X'} \simeq (R^1f'_*K_{X'/C'})^*$ imply $(\det E)^4 = (f'_*\mu^*(f^*K_C \otimes \O_X(D)))^4$. Since $\mu^*(f^*K_C \otimes \O_X(D))$ comes from $C'$, the pull--back gives $({f'}^*\det E)^4 \simeq (\mu^*(f^*K_C \otimes \O_X(D)))^4$. Using \Formel{diag}, we find that $\mu^*(K_X \otimes (f^*K_C \otimes \O_X(D))^{-2})$ is a torsion line bundle on $X'$. Hence
 \[T = K_X \otimes (f^*K_C \otimes \O_X(D))^{-2}\]
is torsion on $X$. After an \'etale cover $\nu: \tilde{X} \to X$, $T$ becomes trivial. If $\tilde{f}: \tilde{X} \to \tilde{C}$ is the Stein factorization of $\nu \circ f$, then $K_{\tilde{X}} \simeq (\tilde{f}^*K_{\tilde{C}} \otimes \O_{\tilde{X}}(\tilde{D}))^2$, where $\tilde{D}$ is defined as in \Abs{ram}. We may hence assume $K_X \simeq (f^*K_C \otimes \O_X(D))^2$. Claim~3 is proved.

\

If $K_X \simeq (f^*K_C \otimes \O_X(D))^2$, then $K_{X'/C'} \simeq \mu^*(f^*K_C \otimes \O_X(D))$, implying $R^2f'_*\mu^*(f^*K_C \otimes \O_X(D)) \simeq R^2f'_*K_{X'/C'} \simeq \O_{C'}$ by duality. Since $\theta_D$ is an isomorphism, the push--forward of the top line in \Formel{diag} gives
 \[0 \lra R^1f'_*\mu^*\Omega_X^1 \lra R^1f'_*\Omega_{X'/C'}^1  \simeq E \otimes E^* \lra \O_{C'} \lra 0.\]
The second map is nothing but the trace map. The pull--back of an ample divisor on $X$ to $X'$ induces a non--trivial section of $R^1f'_*\mu^*\Omega_X^1$. We get a non--nilpotent endomorphism of $E$ whose trace is zero. This shows that $E$ is decomposable (\cite{Ati}, Proposition~16). Claim~1 implies that $E$ decomposes into a sum of line bundles of the same degree. Since $c_1(\Omega_{X'/C'}^1) = \frac{1}{2}c_1(\mu^*K_X) = c_1(\mu^*(f^*K_C \otimes \O_X(D))) = c_1({f'}^*(\sigma^*K_C \otimes \O_{C'}(R_D)))$ on $X'$ by Claim~1, we have $\deg(E) = \deg(\sigma^*K_C \otimes \O_{C'}(R_D))$. In particular we have that
 \begin{equation} \label{even}
   \deg(R_D) = m \cdot \mbox{$\sum\nolimits_{i=1}^l \frac{m_i-1}{m_i}$} \quad \mbox{is even.} 
 \end{equation}

\vspace{0.2cm}

\noindent {\bf Claim~4:} {\em There exists an \'etale cover $\tilde{X} \to X$ such that $\tilde{X}$ is a smooth abelian fibration over a smooth curve $\tilde{C}$.} 

\vspace{0.2cm}

\noindent {\em Proof of Claim~4.} By Claim~3, we may assume $K_X \simeq (f^*K_C \otimes \O_X(D))^2$. Then $K_F \simeq N_{F/X}^*$ on the reduction $F$ of any fiber implying that the only singular fibers of $f$ are multiples of hyperelliptic surfaces of multiplicities $2,3,4,6$. Let $s_2, s_3, s_4, s_6$ be their number, respectively.

We call a ramified base--change $\tilde{C} \to C$ {\em \'etale on $X$}, if the induced map $\tilde{X} = (X \times_{\tilde{C}} C)^{norm} \to X$ is \'etale. Assume that we have some even number of points on $C$, over which we have only multiple fibers of even multiplicity $2,4$ or $6$. Let $\tilde{C} \to C$ be the 2:1--cover, ramified along this points. The corresponding base--change is \'etale on $X$ and reduces the multiplicity of the chosen fibers from $2$, $4$, $6$ to $1$, $2$, $3$, respectively. After a base--change, \'etale on $X$, we may hence assume $s_2 + s_4 + s_6 = 0$ or $1$. Similarly, using 3:1--covers, we may assume $s_3 + s_6 = 0,1$ or $2$. 

If $C \not= \PN{1}$, then $C$ admits \'etale covers of arbitrary degree and we may assume that each $s_i$ is divisible by any given number. Hence, in this case, a ramified base--change, {\'e}tale on $X$, leads to a smooth abelian fibration $\tilde{f}: \tilde{X} \to \tilde{C}$.

If $C = \PN{1}$, then the total number of multiple fibers is at least $3$ since $K_X$ is nef (or see \cite{ViZuiso}). As explained above, we may assume $s_2 + s_4 + s_6 = 0, 1$ and $s_3 + s_6 = 0,1, 2$. Then $s_2 + s_4 + s_6 = 1$ and $s_3 + s_6 = 2$, i.e., $s_6 = 0$, $s_3 = 2$ and either $s_2 = 1$ or $s_4 = 1$. In the case $s_3 = 2$, $s_2 = 1$ we have $\deg(R_D) = 11$. In the case $s_3 = 2$, $s_4 = 1$ we have $\deg(R_D) = 25$. By \Formel{even}, on the other hand, $\deg(R_D)$ must be an even number. Claim~4 is proved. 

\vspace{0.2cm}

\noindent Claim~4 together with Claim~1 prove the \Prop{kod1}.
\end{proof2}

%%%%%%%%%%%%%%%%%%

\

\noindent {\sc Universit\"at Bayreuth, Mathematisches Institut}

\noindent {\sc D--95440 Bayreuth, Germany}

\noindent {\em e--mail: 

ivo.radloff@uni-bayreuth.de

priska.jahnke@uni-bayreuth.de}

\end{document}